\documentclass{amsart}
\usepackage{amsmath, amsthm, latexsym, amssymb, amsfonts}\usepackage[pdftex]{graphicx,color} 

\usepackage[colorlinks, bookmarksdepth=3,citecolor=blue]{hyperref}

\usepackage{tikz}
\usetikzlibrary{patterns}
\usetikzlibrary{decorations.markings}
\tikzset{->-/.style={decoration={
  markings,
  mark=at position #1 with {\arrow{>}}},postaction={decorate}}}

\newenvironment{pf}{\proof[\proofname]}{\endproof}
\theoremstyle{plain}
\newtheorem{Th}{Theorem}[section]
\newtheorem{Conj}[Th]{Conjecture}

\newtheorem{Prop}[Th]{Proposition}

\numberwithin{equation}{section}
\numberwithin{figure}{section}
\theoremstyle{definition}

\newtheorem{Ex}[Th]{Example}

\newtheorem{Prob}[Th]{Problem}

\newcommand{\cal}[1]{\mathcal{#1}}

\newcommand{\C}{\mathbb C}

\newcommand{\N}{\mathbb N}
\newcommand{\Z}{\mathbb Z}
\newcommand{\R}{\mathbb R}

\newcommand{\cO}{\cal O}

\newcommand{\cR}{\cal R}

\newcommand{\la}{\langle}
\newcommand{\ra}{\rangle}

\newcommand{\Sig}{\Sigma}

\newcommand{\Img}{\operatorname{Im}}

\newcommand{\conv}{\operatorname{conv}}

\newcommand{\spn}{\operatorname{span}}

\newcommand{\vol}{\operatorname{vol}}
\newcommand{\Vol}{\operatorname{Vol}}
\newcommand{\V}{\operatorname{V}}
\newcommand{\mv}{\operatorname{v}}

\newcommand{\GL}{\operatorname{GL}}

\newcommand{\rs}[1]{Section~\ref{S:#1}}

\newcommand{\rp}[1]{Proposition~\ref{P:#1}}

\newcommand{\re}[1]{(\ref{e:#1})}

\newcommand{\rt}[1] {Theorem~\ref{T:#1}}

\newcommand{\rf}[1]{Figure~\ref{F:#1}}
\newcommand{\rpr}[1]{Problem~\ref{Pr:#1}}

\begin{document}


\title{The volume polynomial of lattice polygons}
\author[Ivan Soprunov]{Ivan Soprunov}
\address[Ivan Soprunov]{Department of Mathematics and Statistics\\ Cleveland State University\\ Cleveland, OH 44115}
\email{i.soprunov@csuohio.edu}
\author[Jenya Soprunova]{Jenya Soprunova}
\address[Jenya Soprunova]{Department of Mathematical Sciences\\ Kent State University\\ Kent, OH 44242}
\email{esopruno@kent.edu}
\keywords{volume polynomial, lattice polytope, mixed volume, integer quadratic form, intersection number, tropical curve, toric surface}
\subjclass[2020]{Primary 52B20, 52A39; Secondary 11H55, 14T10, 14M25}

\maketitle

{\centering\footnotesize{\it Dedicated to Frank Sottile, our mentor and friend, on the occasion of his 60th birthday.\par}}

\begin{abstract} 
We prove that every indefinite quadratic form with non-negative integer coefficients is the volume polynomial of a pair of lattice polygons. This solves the discrete version of the Heine--Shephard problem for two bodies in the plane. As an application, we
show how to construct a pair of planar tropical curves (or a pair of divisors on a toric surface) with given intersection number and self-intersection numbers.
\end{abstract}


\section{Introduction}

With any collection $K_1,\dots, K_n$ of convex bodies in $\R^d$ and non-negative scalars $x_1,\dots, x_n\in\R_{\geq 0}$
one can associate a convex body
$$x_1K_1+\dots+x_nK_n=\{x_1a_1+\dots+x_na_n : a_i\in K_i, 1\leq i\leq n\}\subset\R^d.$$
In \cite{Min03} Minkowski showed that the $d$-dimensional volume of this body depends polynomially on the scalars, that is, 
$\vol_d(x_1K_1+\dots+x_nK_n)$ is a homogeneous  degree $d$ polynomial in $x_1,\dots,x_n$. 
This polynomial is called the {\it volume polynomial} of $K_1,\dots, K_n$ and its coefficients are the mixed volumes 
(up to multinomial coefficients).  In the case of two planar convex bodies $K,L$  the volume polynomial is a quadratic form
\begin{equation}\label{e:vol2}
\vol_2(xK+yL)=\vol_2(K)x^2+2\mv(K,L)xy+\vol_2(L)y^2,
\end{equation}
 where the middle coefficient $\mv(K,L)$ is called the {\it mixed volume} of $K$ and $L$. It can be expressed by evaluating  \re{vol2} 
 at $x=y=1$ as
  \begin{equation}\label{e:mixed}
  \mv(K,L)=\frac{1}{2}\left(\vol_2(K+L)-\vol_2(K)-\vol_2(L)\right).
  \end{equation}

The classical Minkowski inequality provides a relation between the coefficients of \re{vol2}:
\begin{equation}\label{e:mink}
\vol_2(K)\vol_2(L)\leq \mv(K,L)^2.
  \end{equation}
In other words, \re{vol2} is an indefinite quadratic form. It is not hard to show that every indefinite quadratic form
with non-negative coefficients is the volume polynomial of some planar convex bodies $K,L$. We give a short
proof of this fact in \rp{real}. It also follows from a more general construction of Shephard \cite{Shephard1960}.

The study of polynomial inequalities between the coefficients of the volume polynomial is the core of the Brunn-Minkowski
theory of convex bodies. We refer to the book of Schneider \cite{Schneider2014} for the most complete account of this theory. Still the problem of giving an explicit description for the space of volume polynomials in terms of coefficient inequalities is wide open.
Besides the case of two planar bodies mentioned above, such a description is only known for three planar bodies
($n=3, d=2$) provided by Heine \cite{heine1938wertvorrat}, and two bodies in any dimension ($n=2, d\in\N$) provided by Shephard \cite{Shephard1960}. We recall their descriptions in \rs{conclusion}.
Inequalities such as the Aleksandrov-Fenchel inequality and Shephard's determinantal inequalities uncover
deep connections between mixed volumes, but yet do not provide a complete description for the space of volume polynomials, in general, see \cite{Shephard1960, ABS20}. Some new inequalities describing the square-free part of the volume polynomial for $n=4$ and $d=2$ have been recently found in \cite{AS23}.

In this note we consider a discrete version of the Heine--Shephard problem. Namely, we are interested in describing the space of
volume polynomials of {\it lattice polytopes}. Recall that a convex polytope $P\subset\R^d$ is a lattice polytope if its vertices belong to the integer lattice $\Z^d$. In this case it is appropriate to normalize the Euclidean volume $\vol_d$ by a factor of $d!$. 
Then the normalized volume and mixed volume take integer values on lattice polytopes. 
We use $\Vol_d$ and $\V$ to denote the normalized volume and mixed volume functions, respectively. 

\begin{Prob}\label{Pr:main}
Let $\Vol_d$ be the normalized volume. Describe the set of the volume polynomials
$\Vol_d(x_1P_1+\dots+x_nP_n)\in\Z[x_1,\dots,x_n]$ over all collections of lattice polytopes $P_1,\dots, P_n$ in $\R^d$ in terms
of polynomial inequalities for the coefficients, if possible.
\end{Prob}

In this paper we provide a solution to \rpr{main} in the smallest non-trivial case $n=d=2$. We prove that
every indefinite quadratic form with non-negative integer coefficients is the volume polynomial of a pair of lattice polytopes in $\R^2$,
see \rt{main}. Our proof is constructive. An implementation in Magma \cite{Magma} can be found at
\url{https://github.com/isoprou/volume_polynomial}.

One can also view \rpr{main} as a discrete Blaschke-Santal\'o problem. Motivated by the work of Blaschke \cite{Blaschke1916}, Santal\'o studied the map which assigns to every planar convex body a triple of its geometric invariants,
such as area, perimeter, diameter, width etc., \cite{Santalo61}. The image of such a map is called a Blaschke-Santal\'o diagram. 
Now let $\phi$ be the map which sends a pair of convex planar bodies $K,L$ to the triple $(\Vol_2(K), \V(K,L),\Vol_2(L))$. 
Then the diagram $\Img\phi$ is the closed semialgebraic set 
\begin{equation}\label{e:diagram}
\Img\phi=\{(x,y,z) \in\R_{\geq 0}^3: y^2\geq xz\}
\end{equation}
as follows from Minkowski inequality \re{mink} and the remark after that. 
Restricting $\phi$ to the set of pairs of lattice polytopes in the plane, we obtain the {\it discrete diagram} which, according to
our result in \rt{main}, is the set of lattice points of the semialgebraic set \re{diagram}. In a similar spirit, Scott \cite{Scott76} described
the discrete diagram for the set of coefficients of the Ehrhart polynomial of lattice polygons. For the actual picture of the diagram see, for example, \cite[Thm. 2.1]{RootsEhrhart}.
 
There is a strong motivation to consider mixed volumes of lattice polytopes coming from the intersection theory in toric and tropical geometry. In \rs{application} we discuss an interpretation of our result in terms of the
 intersection numbers of planar tropical curves and divisors on toric surfaces. Finally, in \rs{conclusion} we 
give a conjectural description of the normalized volume polynomial for three lattice polygons and for two polytopes in higher dimensions,
based on Heine's and Shephard's results for convex bodies.
 
\subsection*{Acknowledgments} 
 We are grateful to Gennadiy Averkov for fruitful discussions and to the anonymous referees for their careful reading of the manuscript and for providing insightful comments. 

\section{Preliminaries}

In this section we recall basic facts about mixed volumes of convex bodies. For the most part we restrict ourselves to
the case of lattice polytopes in $\R^2$, as this is the focus of this paper. For the general theory of mixed volumes we refer to \cite[Sec. 5.1]{Schneider2014}. 

A {\it convex body} is a non-empty compact convex set $K\subset\R^d$. The {\it support function} of $K$,
$h_K:\R^d\to\R$ is defined by $h_K(u)=\max\{\la u,v\ra : v\in K\}$, where $\la u,v\ra$ is the standard inner product in $\R^d$. A {\it polytope} $P\subset\R^d$ is the convex hull of finitely many points in $\R^d$. A {\it lattice polytope} $P\subset\R^d$ is the convex hull of finitely many points in $\Z^d$. 
The {\it surface area measure} of $P$ is the discrete measure supported on the set of outer unit normals to 
the facets of $P$ whose values are the $(d-1)$-dimensional volumes of the corresponding facets.

In general, the mixed volume is a polarization of the usual volume, that is, it
is the unique symmetric and multilinear function of $d$ convex bodies which coincides with the volume when the $d$ bodies are the same \cite[Sec. 5.1]{Schneider2014}. In the case of dimension $d=2$, the polarization formula is given in \re{mixed}.
Let $P,Q$ be polytopes in $\R^2$. As seen from \re{mixed}, the mixed volume $\V(P,Q)$ is invariant under simultaneous linear transformations with determinant $\pm 1$ and under independent translations of $P$ and $Q$. 
There is a formula for $\V(P,Q)$ in terms of the support function of $P$ and the surface area measure of $Q$ (see \cite[Thm. 5.1.7]{Schneider2014}) which, in particular, implies the non-negativity of the mixed volume.
Below we present a lattice version of this formula which we will use in
\rs{main}.

Assume $P,Q\subset\R^2$ are lattice polytopes. As mentioned in the introduction, for lattice polytopes we normalize the area by a factor of 2. For example, a lattice rectangle $P=[0,a]\times [0,b]$ for $a,b\in\Z_{\geq 0}$
has normalized volume $\Vol_2(P)=2ab$ and a right triangle $Q$ with vertices $(0,0)$, $(a,0)$, and $(0,b)$ has
$\Vol_2(Q)=ab$. Furthermore, if $I\subset\R^2$ is a lattice segment then we
define its normalized 1-dimensional volume as its lattice length, i.e. $\Vol_1(I)=|I\cap\Z^2|-1$. 

When defining the surface area measure for lattice polytopes it is more convenient to work with primitive normals
rather than unit normals. A vector $u\in\Z^2$ is {\it primitive} if its coordinates are relatively prime.
Let $U_P$ be the set of outer primitive normals to the sides of $P$ and $S_P$ the corresponding surface
area measure supported on $U_P$ with values $S_P(u)=\Vol_1(P^{u})$. Here $P^u$ represents 
the side of $P$ with outer primitive normal $u$. Then the normalized mixed volume can be computed as follows
(see \cite[Sec. 2.5]{mvClass2019} for details)

\begin{equation}\label{e:formula-int}
\V(P,Q)=\sum_{u\in U_Q} h_P(u)\Vol_1(Q^u).
\end{equation}

Note that $h_P(u)\Vol_1(Q^u)$ equals the inner product of a vertex $v\in P$
and a scalar multiple of the primitive vector $u$. We give a small example below.

\begin{Ex} Let $P=\conv\{(0,0),\,(0,2),\,(4,1),\,(4,0)\}$ and $Q=\conv\{(0,0),\,(0,2),\,(3,0)\}$. 
Then $U_Q$ consists of three primitive vectors: $(2,3)$, $(0,-1)$, and $(-1,0)$.
In \rf{example} we show $Q$ with the primitive vectors rescaled
by the lattice lengths of the corresponding sides, $u_1=(2,3)$, $u_2=(0,-3)$, and $u_3=(-2,0)$.
 \begin{figure}[h]
\begin{center}
\includegraphics[scale=.32]{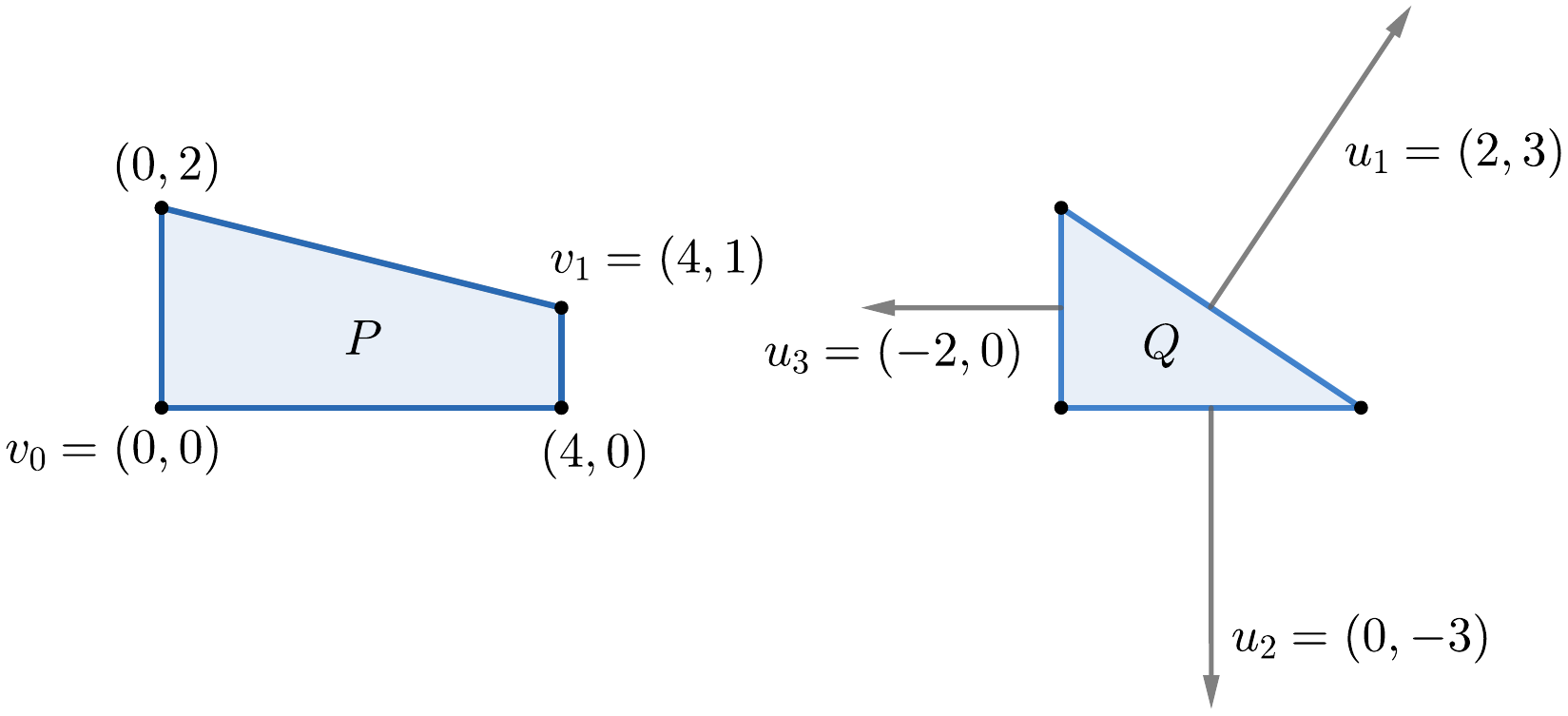} 
\caption{The vertices of $P$ and the surface area measure of $Q$}
\label{F:example} 
\end{center}
\end{figure}
The inner product $\la u_i,v\ra$ attains its maximum at the vertex $v_1=(4,1)\in P$ for $i=1$ and 
at $v_0=(0,0)\in P$ for $i=2,3$. Therefore, by \re{formula-int},
$$\V(P,Q)=\la u_1,v_1\ra + \la u_2,v_0\ra +\la u_3,v_0\ra=11+0+0=11.$$
\end{Ex}

\section{Quadratic forms}

In this section we review basic notions in the theory of binary quadratic forms over the integers. We then prove a 
reduction algorithm for indefinite forms similar to the classical reduction for positive definite forms (see, for example, \cite[Ch. 2]{Buell}).

Consider a binary quadratic form
$$f(x,y)=ax^2+2bxy+cy^2,
$$
where $a,b,c$ are integers. It is called {\it indefinite} if its discriminant $\Delta=4(b^2-ac)$ is non-negative. We
can express $f$ in a matrix form
$$f(x,y)=\left(\begin{matrix}x&y\end{matrix}\right)M\left(\begin{matrix}x \\ y\end{matrix}\right),$$
where $M=\left(\begin{matrix}a & b \\ b & c\end{matrix}\right)$ is an integer symmetric matrix. 

Recall that the unimodular group $\GL(2,\Z)$ is the group of invertible 
integer matrices $G$ with $\det G=\pm 1$. Two forms $f$ and $f'$ are called {\it equivalent} if they are related by a unimodular change of variables.
 In other words, $f'$ is equivalent to $f$ whenever
$$f'(x,y)=\left(\begin{matrix}x&y\end{matrix}\right)G^TMG\left(\begin{matrix}x \\ y\end{matrix}\right),$$
for some $G\in\GL(2,\Z)$. Here $G^T$ is the transpose of $G$. In this case we say that $G$ transforms
$f$ to $f'$.
More explicitly, if 
$G=\left(\begin{matrix}x_1 & x_2 \\ y_1 & y_2\end{matrix}\right)$ then
$f'(x,y)=a'x^2+2b'xy+c'y^2$ with
\begin{align*}
a'&=ax_1^2+2bx_1y_1+cy_1^2=f(x_1,y_1)\\
b'&=ax_1x_2+b(x_1y_2+x_2y_1)+cy_1y_2\\
c'&=ax_2^2+2bx_2y_2+cy_2^2=f(x_2,y_2).
\end{align*}
Clearly, equivalent forms have the same discriminant.

Our next result is a key tool for constructing lattice polytopes with a given volume polynomial in \rs{main}.

\begin{Prop}\label{P:reduction} Every integer indefinite binary form with positive coefficients is equivalent to a form $f(x,y)=ax^2+2bxy-cy^2$ for some $a,b,c\in\Z_{\geq 0}$ such that $f(1,1)=a+2b-c>0$.
Moreover, there exists a matrix
$G=\left(\begin{matrix}x_1 & x_2 \\ y_1 & y_2\end{matrix}\right)\in\GL(2,\Z)$ satisfying $x_i\geq y_i\geq 0$ for $i=1,2$
which transforms $f$ to the original form.
\end{Prop}

\begin{pf}
Consider $F(x,y)=Ax^2+2Bxy+Cy^2$ with $A,B,C\in\N$ and let $k=B^2-AC\geq 0$.
If $A<C$ then we apply the matrix $\left(\begin{matrix}0 & 1 \\ 1 & 0\end{matrix}\right)$ which 
swaps the variables. Thus, we may assume $A\geq C$. Since $k\geq 0$ this implies that $B\geq C$.

Next,  we divide $B$ by $C$ with a remainder, $B=sC+r$ where $s\geq 1$ and $0\leq r<C$. Applying the matrix
$\left(\begin{matrix}0 & 1 \\ 1 & -s\end{matrix}\right)$ we obtain a form $ax^2+2bxy+cy^2=F(y,x-sy)$ with
$a=C$, $b=r$, and $c=F(1,-s)$. Note that we have $a>b>c$, since $a>b$ and the discriminant $k=b^2-ac\geq 0$. 
While $c>0$ we keep dividing $b$ by $c$ with a remainder. Since throughout this process $c$ is getting strictly smaller,
we will eventually end up with the form $f(x,y)=ax^2+2bxy+cy^2$ where $a>b\geq 0\geq c$. 

At this point it may happen that $f(1,1)\leq 0$. If this is the case, we modify the last step as follows. Let $F(x,y)=Ax^2+2Bxy+Cy^2$ with $A,B,C\in\N$ be the form in the penultimate step and $B=sC+r$ be the last division with remainder, so $F(1,-s)\leq 0$.
Choose $s'$ to be the smallest integer such that $F(1,-s')\leq 0$. Note that $s'\geq 1$ since $F$ has positive coefficients.
Then, applying the matrix $\left(\begin{matrix}0 & 1 \\ 1 & -s'\end{matrix}\right)$ we obtain a form $f(x,y)=ax^2+2bxy+cy^2=F(y,x-s'y)$ with $a=C$, $b=B-s'C$, and $c=F(1,-s')$. Clearly, $a>0$, $c\leq 0$, and, since $s'\leq s$, we also have $b\geq r\geq 0$. Moreover, $f(1,1)=F(1,1-s')>0$, since otherwise $s'$ was not the smallest. Finally, by switching the sign of $c$, we may write the resulting form as
$f(x,y)=ax^2+2bxy-cy^2$ with $a,b,c\in\Z_{\geq 0}$.


It remains to look at the matrix $G$ that transforms $f$ to $F$. It is the product of the inverses of the matrices used above, i.e.
$$G=
\left(\begin{matrix}s_n & 1 \\ 1 & 0\end{matrix}\right)\cdots 
\left(\begin{matrix}s_1 & 1 \\ 1 & 0\end{matrix}\right)
\left(\begin{matrix}0 & 1 \\ 1 & 0\end{matrix}\right),
$$
where $n\geq 1$, $s_i\geq 1$ for $1\leq i\leq n$, and the rightmost matrix may or may not be present. It is easy to see by induction that 
$G$ satisfies the condition in the statement of the proposition.
%
\end{pf}

\section{The volume polynomial}\label{S:main}

In this section we prove our main result (\rt{main}) which describes all normalized volume polynomials of pairs of lattice polytopes in $\R^2$. First, we look at the easier case of the usual volume polynomial of two planar bodies.
Although this is not new and  was probably already known to Minkowski, we include it here for comparison of the complexity of the continuous and the discrete cases. In \cite{Shephard1960} Shephard proved a generalization of this result to arbitrary dimension using a pair of simplices. Our construction uses a pair of coordinate rectangles instead.

\begin{Prop}\label{P:real}
Let $F(x,y)=Ax^2+2Bxy+Cy^2$ be an indefinite quadratic form with
$A,B,C\in\R_{\geq 0}$. Then there exist convex bodies $K,L\subset\R^2$ such that $F(x,y)$ is the volume polynomial
of $K$ and $L$, that is, $F(x,y)=\vol_2(xK+yL)$ for $x,y\geq 0$.
\end{Prop}

\begin{pf} We let $K$ and $L$ be rectangular boxes $K=[0,\alpha]\times[0,\beta]$ and $L=[0,\gamma]\times[0,\delta]$
for some $\alpha, \beta, \gamma, \delta\in\R_{\geq 0}$. Then
$\vol_2(K)=\alpha\beta$, $\mv(K,L)=(\alpha\delta+\beta\gamma)/2$, and  $\vol_2(L)=\gamma\delta$.
We will show that the system
$$A=\alpha\beta,\quad B=(\alpha\delta+\beta\gamma)/2,\quad C=\gamma\delta$$
has a non-negative solution. Indeed, first assume $A>0$. Put $\alpha=1$ and so $\beta=A$. Multiplying the middle equation by $\gamma$ and using the last equation we obtain:
$$A\gamma^2-2B\gamma+C=0.$$
Since $k=B^2-AC\geq 0$ this equation has a non-negative solution $\gamma=\frac{B+\sqrt{k}}{A}$. Note that $\gamma=0$ only if $B=0$ and $C=0$. In this case $\delta=0$. Otherwise, $\delta=C/\gamma$. 

Now, if $A=0$ and $B>0$ then we put $\alpha=1$, $\beta=0$, $\gamma=\frac{C}{2B}$, and $\delta=2B$. Finally, if $A=0$ and $B=0$ then
we put $\alpha=0$, $\beta=0$, $\gamma=1$, and $\delta=C$.
\end{pf}

\begin{Th}\label{T:main} Let $F(x,y)=Ax^2+2Bxy+Cy^2$ be an indefinite quadratic form with
$A,B,C\in\Z_{\geq 0}$. Then there exist lattice polytopes $P,Q\subset\R^2$ such that $F(x,y)$ is the normalized volume polynomial
of $P$ and $Q$, that is $F(x,y)=\Vol_2(xP+yQ)$ for $x,y\geq 0$.
\end{Th}

\begin{pf} Suppose $A$ and $C$ are positive, and, hence, so is $B$. Then,
according to \rp{reduction}, there exist a form $f(x,y)=ax^2+2bxy-cy^2$ with $a,b,c\in\Z_{\geq 0}$, and $f(1,1)=a+2b-c>0$,
and a unimodular matrix
$G=\left(\begin{matrix}x_1 & x_2 \\ y_1 & y_2\end{matrix}\right)\in\GL(2,\Z)$ with $x_i\geq y_i\geq 0$, $i=1,2$
which transforms $f$ to $F$. Then
$$A=f(x_1,y_1),\quad B=ax_1x_2+b(x_1y_2+x_2y_1)-cy_1y_2,\quad C=f(x_2,y_2).$$ 
Below we show how to construct lattice polytopes $P$ and $Q$
satisfying $\Vol_2(P)=A$, $\V(P,Q)=B$, and $\Vol_2(Q)=C$. The construction differs slightly depending on the parity
of $a$ and $c$.

\noindent {\it Case 1.} Suppose first that $a$ and $c$ are even, so that $a=2d$ and $c=2e$ for some $d,e\in\Z_{\geq 0}$. Then we can write $f(x,y)=ax^2+2bxy-cy^2=2x(dx+by)-2ey^2$.
We view this expression as the normalized area of a rectangle with sides $x$ and $dx+by$ with two opposite corners cut off; each
corner being a right triangle with legs $y$ and $ey$. Thus, we define 
$$P=\conv\{\big(0,y_1\big),\, \big(0,x_1\big),\, \big(ey_1,0\big),\, \big(dx_1+by_1, 0\big),\, \big(dx_1+by_1,x_1-y_1\big),\, \big(dx_1+(b-e)y_1,x_1\big)\}.$$ 
We define $Q$  in the same way replacing  $x_1$ with $x_2$ and  $y_1$ with $y_2$ above.
Since $x_i\geq y_i\geq 0$ and 
$$
dx_i+(b-e)y_i\geq y_i(d+b-e)=y_i(a+2b-c)/2\geq 0,
$$ 
the relative positions of the vertices of $P$ and $Q$ are as shown in Figure~\ref{F:mixed1}.   
Clearly, $\Vol_2(P)=2(dx_1+by_1)-2ey_1^2=f(x_1, y_1)=A$ and, similarly, $\Vol_2(Q)=f(x_2,y_2)=C$.
In the second diagram in Figure~\ref{F:mixed1} we depict the surface area measure of $Q$. This allows us to compute the mixed volume using \re{formula-int}.
 \begin{figure}[h]
\begin{center}
\includegraphics[scale=.34]{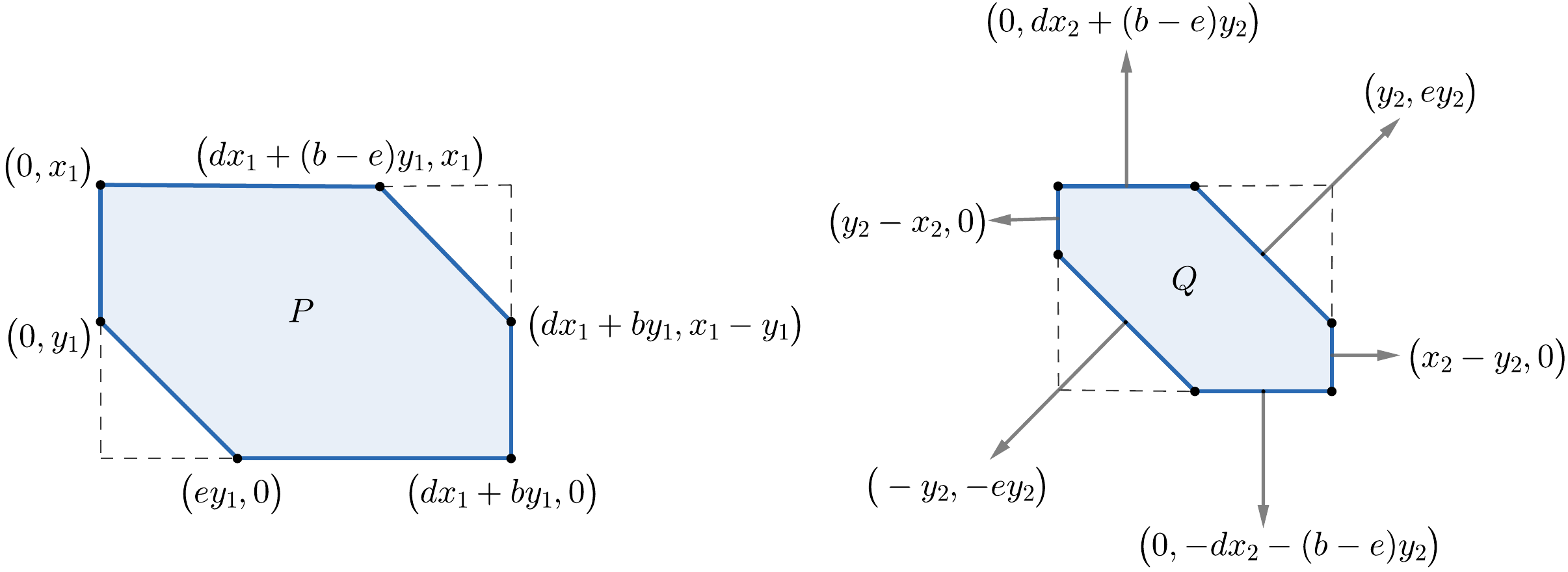} 
\caption{Case 1: $a=2d$ and $c=2e$}
\label{F:mixed1} 
\end{center}
\end{figure}
\begin{align*}
\V(P,Q)&=(x_2-y_2)(dx_1+by_1)+y_2(dx_1+by_1)+ey_2(x_1-y_1)+x_1(dx_2+(b-e)y_2)\\
&-ey_1y_2=ax_1x_2+b(x_1y_2+x_2y_1)-cy_1y_2=B.
\end{align*} 

\noindent {\it Case 2.} If $a=2d$ and $c=2e+1$ for $d,e\in\Z_{\geq 0}$, we construct $P$ and $Q$ as above except the primitive outer normal 
$(1,e)$ is now replaced with $(1,e+1)$. Then
$$P=\conv\{\big(0,y_1\big),\, \big(0,x_1\big),\, \big(ey_1,0\big),\, \big(dx_1+by_1, 0\big),\, \big(dx_1+by_1,x_1-y_1\big),\, \big(dx_1+(b-e-1)y_1,x_1\big)\},$$ 
and $Q$ is obtained by replacing the indices. We have $x_i\geq y_i\geq 0$ and 
$$ dx_i+(b-e-1)y_i\geq (d+b-e-1)y_i=y_i(a+2b-c-1)/2\geq 0,$$ and we get a diagram depicted in Figure~\ref{F:mixed2}.
\begin{figure}[h]
\begin{center}
\includegraphics[scale=.34]{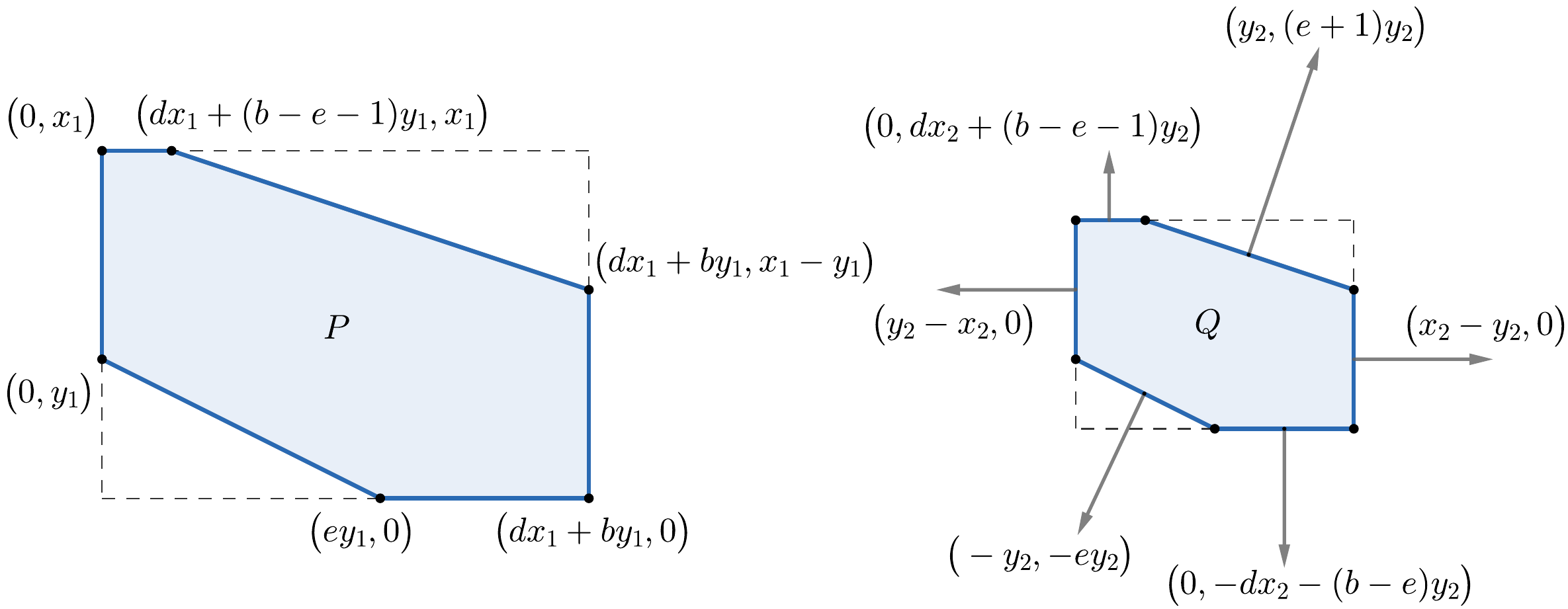} 
\caption{Case 2: $a=2d$ and $c=2e+1$}
\label{F:mixed2} 
\end{center}
\end{figure}
As before, we have $\Vol_2(P)=A, \Vol_2(Q)=C$, and
\begin{align*}
\V(P,Q)&=(x_2-y_2)(dx_1+by_1)+y_2(dx_1+by_1)+(e+1)y_2(x_1-y_1)\\
&+x_1(dx_2+(b-e-1)y_2)-ey_1y_2=ax_1x_2+b(x_1y_2+x_2y_1)-cy_1y_2=B.
\end{align*} 

\noindent {\it Case 3.} Suppose $a=2d+1$ and $c=2e$ for $d, e\in\Z_{\geq 0}$. We now write 
$$f(x,y)=ax^2+2bxy-cy^2=2x(dx+by)+x^2-2ey^2.$$
This represents the normalized area of a trapezoid, which is the union of an $x$ by $dx+by$ rectangle and
an isosceles right triangle with leg $x$, with two opposite corners cut off, where
the corners have normalized area $ey^2$ each. In other words, we let $P$ have the vertices
$$
(0,y_1),\, (0,x_1),\, (ey_1,0),\, ((d+1)x_1+by_1, 0),\, (dx_1+(b+1)y_1,x_1-y_1),\,
(dx_1+(b-e)y_1,x_1),
$$
and $Q$ is defined accordingly. 
Similarly to the above, $x_i\geq y_i\geq 0$ and 
$$dx_i+(b-e)y_i\geq y_i(d+b-e)=y_i(a+2b-c-1)/2\geq 0,$$ and we get the diagram  in Figure~\ref{F:mixed3}.
\begin{figure}[h]
\begin{center}
\includegraphics[scale=.34]{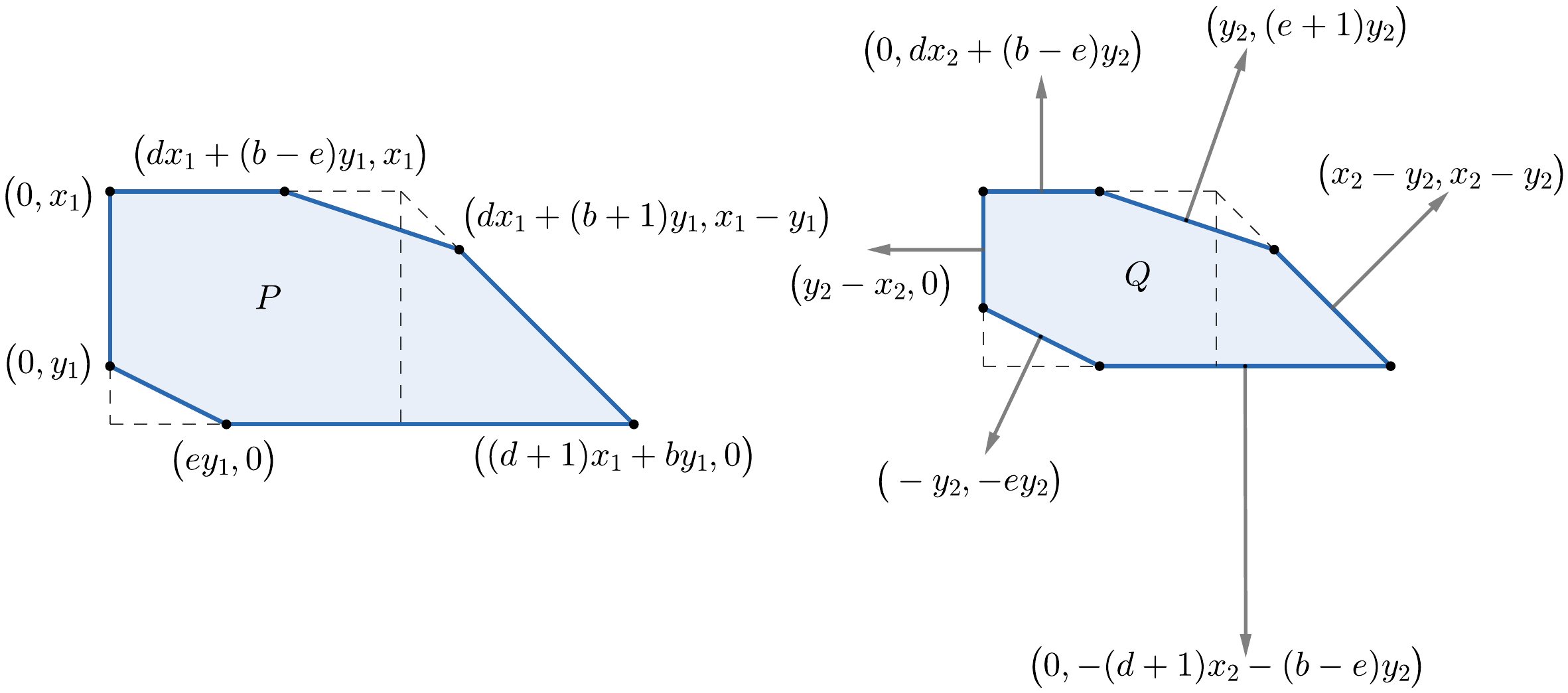} 
\caption{Case 3: $a=2d+1$ and $c=2e$}
\label{F:mixed3} 
\end{center}
\end{figure}
We have $\Vol_2(P)=A$, $\Vol_2(Q)=B$, and 
\begin{align*}
\V(P,Q)&=(x_2-y_2)((d+1)x_1+by_1)+y_2(dx_1+(b-e)y_1)+(e+1)x_1y_2\\
&+x_1(dx_2+(b-e)y_2)-ey_1y_2=ax_1x_2+b(x_1y_2+x_2y_1)-cy_1y_2=B.
\end{align*} 

\noindent {\it Case 4.} Finally, suppose $a=2d+1$ and $c=2e+1$ for $d,e\in\Z_{\geq 0}$. Then $x_i\geq y_i\geq 0$ and 
$$dx_i+(b-e)y_i\geq y_i(d+b-e)=y_i(a+2b-c)/2\geq 0,$$ 
and  we define  each of $P$ and $Q$ as the convex hull of 
$$
 \big(0,y_i\big),\, \big(0,x_i\big),\, \big(ey_i,0\big),\, \big((d+1)x_i+by_i, 0\big),\, \big(dx_i+(b+1)y_i,x_i-y_i\big),\,
\big(dx_i+(b-e-1)y_i,x_i\big),
$$
as depicted in Figure~\ref{F:mixed4}.
\begin{figure}[h]
\begin{center}
\includegraphics[scale=.34]{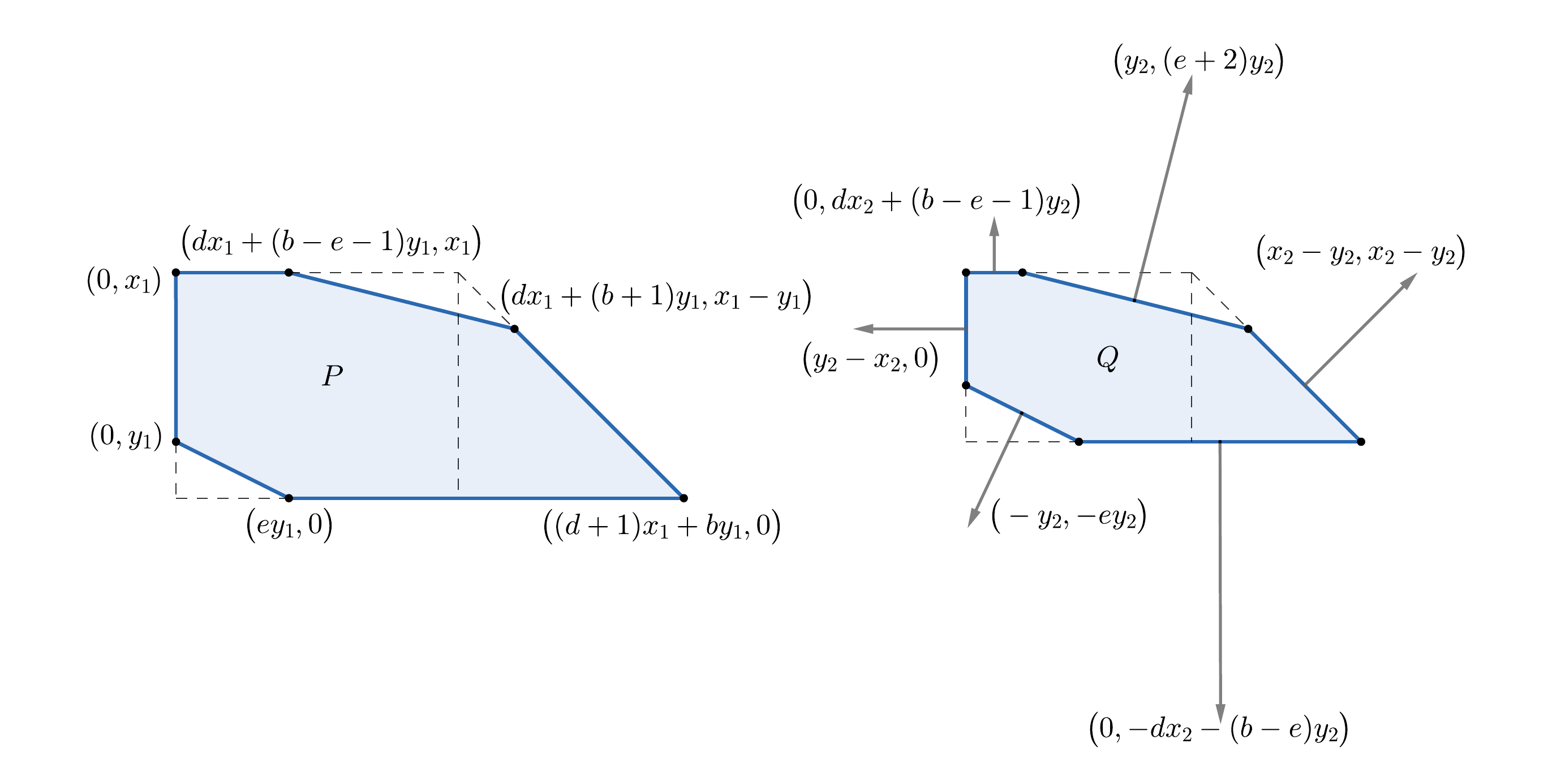} 
\caption{Case 4: $a=2d+1$ and $c=2e+1$}
\label{F:mixed4} 
\end{center}
\end{figure}
Then $\Vol_2(P)=A$, $\Vol_2(Q)=C$, and 
\begin{align*}
\V(P,Q)&=(x_2-y_2)((d+1)x_1+by_1)+y_2(dx_1+(b-e-1)y_1)+(e+2)x_1y_2\\
&+x_1(dx_2+(b-e-1)y_2)-ey_1y_2=ax_1x_2+b(x_1y_2+x_2y_1)-cy_1y_2=B.
\end{align*} 

In remains to consider the case when $A$ or $C$ (or both) is zero. Without loss of generality, assume $A=0$.
Then we take $P$ to be the unit segment $P=\conv\{(0,0),\, (1,0)\}$. To construct $Q$ we write $C=sB-r$, for
some $s,r\in\Z$ with $s\geq 1$ and $0< r\leq B$, and let
$$Q=\conv\{(s,0),\,(1,0),\,(0,r),\,(0,B)\}.$$
Then $\Vol_2(Q)=C$ and $\V(P,Q)=B$, as required. Note that when, in addition, $C=0$ we have $s=1$, $r=B$ and, hence, $Q$ is also a segment.  

\end{pf}

\section{Connection to intersection numbers of curves}\label{S:application}

In this section we discuss an implication of \rt{main} for intersection numbers of plane tropical curves, as well as intersection numbers of
divisors on toric surfaces. We recall basic facts below and refer to \cite{Fulton, TorVar, Tropical} for the general theory of toric and tropical varieties.

\subsection{Planar tropical curves} A tropical bivariate polynomial $f$ is a piece-wise linear function 
$f(x,y) = \min\{ c_{a,b} + a x + b y  : (a,b) \in S\}$ where $S\subset\Z^2$ is a finite set, called the support of $f$, and $c_{a,b}\in\R$.
The convex hull $P_f=\conv S$ is called the {\it Newton polytope} of $f$.

The set of points $(x,y)\in\R^2$ where $f$ is not differentiable (that is where the minimum in the definition of $f$ is attained at least twice)
is called the (planar) {\it tropical curve} $V_f \subset \R^2$ associated to $f$. Geometrically, $V_f$ is a polyhedral complex of pure dimension 1. There is a duality between the polyhedral complex $V_f$ and
the regular subdivision $\cR_f$ of $P_f$ induced by lifting $(a,b)\in S$ to height $c_{a,b}$. The vertices (0-dimensional cells) of $V_f$ correspond to the polygons (2-dimensional cells) of $\cR_f$ and the edges (1-dimensional cells) of $V_f$ correspond to the edges
of $\cR_f$. Given an edge $e\in V_f$, its {\it weight} $w(e)$ is the lattice length of the corresponding edge in~$\cR_f$.

Two tropical curves $V_f$ and $V_g$ are said to intersect {\it transversally} if $V_f \cap V_g$ is finite and every point $p \in V_f \cap V_g$ lies in the relative interior of some edge $e$ of $V_f$ and some edge $h$ of $V_g$. In this case, define the {\it local intersection number} $(V_f,V_g)_p$ at point $p \in V_f \cap V_g$
to be $w(e)w(h) |\det(u,v)|$, where $u$ and $v$ are primitive vectors parallel to $e$ and $h$, respectively. Then the {\it intersection number} 
is defined as the sum of local intersection numbers
$$(V_f,V_g)=\sum_{p\in V_f \cap V_g}(V_f,V_g)_p.$$

The following theorem is known as the tropical Bernstein-Khovanskii-Kushnirenko theorem in dimension two, \cite[Thm. 4.6.8]{Tropical}.

\begin{Th}\label{T:trop-BKK} 
The intersection number of two planar tropical curves $V_f$ and $V_g$  intersecting transversally equals the normalized
mixed volume of their Newton polytopes $\V(P_f,P_g)$.
\end{Th}

Let $V_f$ be a tropical curve and let $V_{f'}$ be its translate by $(x_0,y_0)\in\R^2$. In other words  $V_{f'}$ is associated to the tropical polynomial  $f'(x,y)=f(x+x_0,y+y_0)$ which has the same support and, hence, the same Newton polytope as $f$. 
Note that the lifting function $c'_{a,b}$  differs from $c_{a,b}$ by a linear function, $c'_{a,b}=c_{a,b}+ax_0+by_0$, hence, the regular subdivisions $\cR_f$ and $\cR_{f'}$ are also the same.  For almost all values of $(x_0,y_0)$ the curves $V_f$  and $V_{f'}$ intersect transversally, in which case we define  $(V_f,V_f):=(V_f,V_{f'})$ and call it the {\it self-intersection number} of $V_f$. By \rt{trop-BKK}, 
$(V_f,V_f)=\V(P_f,P_f)=\Vol_2(P_f)$. In particular, it does not depend on the choice of $V_{f'}$.

The following is a direct consequence of Theorems \ref{T:main} and \ref{T:trop-BKK}.

\begin{Th}\label{T:trop-main} Given $(A,B,C)\in\Z_{\geq 0}^3$ satisfying $AC\leq B^2$, there exist planar tropical curves $V_f$ and $V_g$ with 
$(V_f,V_f)=A$, $(V_f,V_g)=B$, and $(V_g,V_g)=C$.  
\end{Th} 

\begin{pf} Let $P$ and $Q$ be the lattice polytopes constructed in \rt{main} from the triple $(A,B,C)$. Consider tropical polynomials
$f$ and $g$ with supports $P\cap\Z^2$ and $Q\cap\Z^2$, respectively, and generic coefficients $c_{a,b}$. Then the corresponding 
tropical curves $V_f$ and $V_g$ will intersect transversally and $(V_f,V_f)=\Vol_2(P)=A$, $(V_f,V_g)=\V(P,Q)=B$, and $(V_g,V_g)=\Vol_2(Q)=C$.
\end{pf}

\subsection{Divisors on toric surfaces} 
A toric surface $X$ is an algebraic surface containing the torus $(\C^*)^2$ as a Zariski open subset whose action on itself extends
to the action on $X$. Basic examples of complete toric surfaces include the projective plane, the product of projective lines, and the Hirzebruch surface. One can describe the affine charts of $X$ via a rational polyhedral fan $\Sig$ in $\R^2$ whose 1-dimensional cones (rays) are generated by a collection of primitive vectors $\{u_1,\dots,u_r\}$. Each ray corresponds to a 1-dimensional orbit in $X$ whose closure defines a  torus invariant prime divisor $D_i$ on $X$. Let $D=\sum_{i=1}^r a_iD_i$ be a torus invariant divisor. It defines a rational polytope\footnote{Since we work with outer normals rather than inner normals, we modify the standard definition accordingly.}
$P_D=\{x\in\R^2 : \la x,u_i\ra\leq a_i, 1\leq i\leq r\}$. 

As shown in \cite[Thm. 4.2.8]{TorVar}, $D=\sum_{i=1}^r a_iD_i$ is a Cartier divisor if and only if there is a continuous piece-wise linear function $\phi_D$ on $\R^2$ which is linear on each cone of $\Sig$ and $\phi_D(u_i)=a_i$ for $1\leq i\leq r$. The global sections of the corresponding line bundle $\cO(D)$ can be identified with  the space of Laurent polynomials supported on $P_D\cap\Z^2$, i.e. $\Gamma(X,\cO(D))\cong \spn_\C\{x_1^{a_1}x_2^{a_2} : (a_1,a_2)\in P_D\cap\Z^2\}$, see \cite[Thm. 4.3.3]{TorVar}.
Furthermore, $\cO(D)$ is globally generated if and only if $\phi_D$ is convex in which case  $P_D$ is a lattice polytope, see \cite[Thm. 6.1.7]{TorVar}.

We next recall a version of the Bernstein-Khovanskii-Kushnirenko theorem (see \cite[Sec. 5.4-5.5]{Fulton}) in the case of intersection numbers of divisors on toric surfaces. Let $D,E$ be globally generated divisors on $X$ and $f\in\Gamma(X,\cO(D))$, $g\in\Gamma(X,\cO(E))$ be generic sections of the corresponding line bundles. 
By above, $f$ and $g$ are Laurent polynomials supported on $P_D\cap \Z^2$ and $P_E\cap\Z^2$, respectively.
The zero loci of $f$ and $g$ define curves $V_f$ and $V_g$  on $X$ having simple isolated intersections in the torus $(\C^*)^2$. The total number of intersection points is independent of the choice of $f$ and $g$ and equals the intersection number $(D,E)$.

\begin{Th}\label{T:toric-BKK} Let $D,E$ be two globally generated divisors on a toric surface $X$. Then
the intersection number $(D,E)$  equals the normalized mixed volume of their polytopes $\V(P_D,P_E)$.
\end{Th}

We can now interpret the result of \rt{main} as follows.

\begin{Th}\label{T:toric-main} Given $(A,B,C)\in\Z_{\geq 0}^3$ satisfying $AC\leq B^2$, there 
exist a toric surface $X$ and globally generated divisors $D$ and $E$ on $X$ such that 
$(D,D)=A$, $(D,E)=B$, and $(E,E)=C$.  
\end{Th} 

\begin{pf} Let $P$ and $Q$ be the lattice polytopes constructed in \rt{main} from the triple $(A,B,C)$. 
The primitive normals of both $P$ and $Q$ belong to a set of six vectors $\{u_1,\dots, u_6\}$. For example,
in Case 1 the vectors are $\{\pm (1,0), \pm (0,1),\pm (1, e)\}$.
Let $\Sig$ be the fan generated by $\{u_1,\dots, u_6\}$ and $X$ the toric surface associated to $\Sig$. 
Define $a_i=h_P(u_i)$ and $b_i=h_Q(u_i)$ for $1\leq i\leq 6$, where $h_P$ and $h_Q$ are the support
functions of $P$ and $Q$, respectively. Then the divisors $D=\sum_{i=1}^6 a_iD_i$ and $E=\sum_{i=1}^6 b_iD_i$ are globally
generated and have polytopes $P$ and $Q$, respectively. By Theorems \ref{T:toric-BKK} and
\ref{T:main} we have $(D,D)=\Vol_2(P)=A$, $(D,E)=\V(P,Q)=B$, and $(E,E)=\Vol_2(Q)=C$.
\end{pf}

We note that in the above proof $X$ is not necessarily smooth as the fan $\Sig$ may not be unimodular.
(A fan $\Sig\subset\R^2$ is unimodular if the primitive generators of every 2-dimensional cone in $\Sig$ form a basis for $\Z^2$.)
By subdividing the 2-dimensional cones of $\Sig$ if necessary one obtains a unimodular fan $\Sig'$ which corresponds to a smooth toric surface $X'$. The globally generated divisors $D$ and $E$ on $X'$ are constructed from the polytopes $P$ and $Q$ in the same way as in the above proof.

\section{Conclusion}\label{S:conclusion}
As mentioned in the introduction, the Heine-Shephard problem of describing the space of volume polynomials of
$n$ convex bodies in $\R^d$, although simple enough for two planar bodies, becomes challenging already for three bodies in the plane  and two bodies in higher dimensions (the results of Heine and Shephard) and is open in all other cases. Even the existence of a semialgebraic description of the space of volume polynomials is an open question, in general. (Some recent results on the semialgebraicity of a
projection of the space of volume polynomials for any number of planar bodies are contained in \cite{AS23}.)
In the discrete version (\rpr{main}), the challenges present themselves already in the case of two lattice polygons $P,Q\subset\R^2$. It is not surprising that the discrete version is harder, but the fact that there are no extra conditions on the coefficients of $\Vol_2(xP+yQ)$ besides being integers, allows us to suggest that the same situation may persist in larger cases as well. 
 
Based on Heine's description of volume polynomials of triples of planar bodies (see \cite{heine1938wertvorrat}), we propose the following conjecture.

\begin{Conj}\label{Conj:Heine} Let $M$ be a $3\times 3$ symmetric matrix with non-negative integer entries and let
$F(x,y,z)=\left(x \,\, y\,\, z\right)M\left(\begin{matrix}x \\ y\\ z\end{matrix}\right)$ be the corresponding ternary form.
Assume $\det(M)\geq 0$ and $\det(M_I)\leq 0$, where $\det(M_I)$ is the diagonal $2\times 2$ minor indexed by $I\subset\{1, 2, 3\}$.
Then there exist lattice polygons $P,Q,R\subset\R^2$ such that $F(x,y,z)=\Vol_2(xP+yQ+zR)$.
\end{Conj}

Shephard describes the volume polynomials for a pair of bodies in $\R^d$ as polynomials whose coefficients (up to binomial coefficients) form a log-concave sequence, see \cite{Shephard1960}. This suggests the following conjecture in the discrete case.
\rt{main} confirms it for $d=2$.

\begin{Conj}\label{Conj:Shephard} Let $F(x,y)=\sum_{i=0}^d{d\choose i}a_ix^iy^{d-i}$ be a degree $d$ form where $a_0,\dots, a_d$
are non-negative integers satisfying $a_{i-k}a_{j+k}\leq a_ia_j$ for all $0\leq i\leq j\leq d$ and $0\leq k\leq\min(i,d-j)$.
Then there exist lattice polytopes $P,Q\subset\R^d$ such that $F(x,y)=\Vol_d(xP+yQ)$.
\end{Conj}

\bibliographystyle{alpha} 
\bibliography{lit.bib} 

\end{document}